\documentclass[11pt]{amsart}

\usepackage{amsmath}
\usepackage{amssymb}
\usepackage{mathrsfs}
\usepackage{comment}
\usepackage{color}
\usepackage[colorlinks,citecolor=blue,urlcolor=black,linkcolor=black]{hyperref}

\newtheorem{theorem}{Theorem}[section]
\newtheorem*{theorem*}{Theorem}
\newtheorem{claim}[theorem]{Claim}

\newtheorem{lemma}[theorem]{Lemma}

\newtheorem{corollary}[theorem]{Corollary}

\theoremstyle{definition}
\newtheorem{definition}[theorem]{Definition}

\newtheorem{question}[theorem]{Question}
\newtheorem*{question*}{Question}
\theoremstyle{remark}

\newcommand{\cf}{{\rm cf}}

\newcount\skewfactor
\def\mathunderaccent#1#2 {\let\theaccent#1\skewfactor#2
\mathpalette\putaccentunder}
\def\putaccentunder#1#2{\oalign{$#1#2$\crcr\hidewidth
\vbox to.2ex{\hbox{$#1\skew\skewfactor\theaccent{}$}\vss}\hidewidth}}
\def\name{\mathunderaccent\tilde-3 }

\def\smallbox#1{\leavevmode\thinspace\hbox{\vrule\vtop{\vbox
   {\hrule\kern1pt\hbox{\vphantom{\tt/}\thinspace{\tt#1}\thinspace}}
   \kern1pt\hrule}\vrule}\thinspace}

\def\qedref#1{$\qed_{\reforiginal{#1}}$}

\title{Superclub, splitting, separating statements}
\author{Shimon Garti}
\address{Einstein Institute of Mathematics,
 The Hebrew University of Jerusalem,
 Jerusalem 91904, Israel}
\email{shimon.garty@mail.huji.ac.il}

\author{Saharon Shelah}
\address{Institute of Mathematics
 The Hebrew University of Jerusalem
 Jerusalem 91904, Israel
 and  Department of Mathematics
 Rutgers University
 New Brunswick, NJ 08854, USA}
\email{shelah@math.huji.ac.il}
\urladdr{http://www.math.rutgers.edu/\char`\~shelah}

\thanks{Research supported by ISF grant no. 1838/19. This is publication 1236 in the second author's list}

\begin{document}
\let\labeloriginal\label
\let\reforiginal\ref
\def\ref#1{\reforiginal{#1}}
\def\label#1{\labeloriginal{#1}}

\begin{abstract}
We prove that superclub implies $\mathfrak{s}=\omega_1$. More generally, if $\kappa$ is weakly compact then superclub implies $\mathfrak{s}_\kappa=\kappa^+$.
Based on this statement, we separate tiltan from superclub at successors of supercompact cardinals.
We also use the Galvin property in order to separate tiltan from superclub at both successors of regular and successors of singular cardinals.
\end{abstract}

\maketitle

\newpage

\section{Introduction}

One of the cornerstones of modern set theory is the diamond principle, discovered by Jensen in \cite{MR309729}.
A sequence of the form $(A_\alpha\mid\alpha\in\omega_1)$ is a diamond sequence iff $A_\alpha\subseteq\alpha$ for every $\alpha\in\omega_1$ and if $A\subseteq\omega_1$ then $S_A=\{\alpha\in\omega_1\mid A\cap\alpha=A_\alpha\}$ is a stationary subset of $\omega_1$.
More generally, if $\kappa=\cf(\kappa)>\aleph_0$ and $S$ is a stationary subset of $\kappa$ then an $S$-diamond sequence is a sequence of sets $(A_\alpha\mid\alpha\in{S})$ such that $A_\alpha\subseteq\alpha$ whenever $\alpha\in{S}$ and if $A\subseteq\kappa$ then $S_A=\{\alpha\in{S}\mid A\cap\alpha=A_\alpha\}$ is a stationary subset of $\kappa$.

One says that diamond holds (resp., diamond holds at $S$) iff there is a diamond sequence (resp., an $S$-diamond sequence).
This statement is denoted by $\Diamond$ (resp., $\Diamond_S$).
Jensen proved that $\Diamond_S$ holds in the constructible universe at every stationary $S\subseteq\kappa=\cf(\kappa)>\aleph_0$.

A related combinatorial principle was discovered by Ostaszewski at nearly the same time.
This is the club principle which is used in \cite{MR438292} to build a certain example of a topological space.
For avoiding ambiguity we shall call this principle \emph{tiltan}.
A tiltan sequence is a sequence of the form $(T_\delta\mid\delta\in\lim(\omega_1))$ such that every $T_\delta$ is an unbounded subset of $\delta$ and for every $A\in[\omega_1]^{\omega_1}$ the set $T_A=\{\delta\in\lim(\omega_1)\mid T_\delta\subseteq A\cap\delta\}$ is a stationary subset of $\omega_1$.
Tiltan holds iff there exists a tiltan sequence, and this statement is denoted by $\clubsuit$.
Once again, if $\kappa=\cf(\kappa)>\aleph_0$ and $S$ is a stationary subset of $\kappa$ consists of limit ordinals, then $\clubsuit_S$ means that there is a tiltan sequence $(A_\delta\mid\delta\in{S})$.
That is to say that for every $A\in[\kappa]^\kappa$ the set $T_A=\{\delta\in{S}\mid A\cap\delta\supseteq A_\delta\}$ is stationary in $\kappa$.

It follows directly from the definitions that $\Diamond_S$ implies $\clubsuit_S$.
However, $\clubsuit_S$ is strictly weaker than $\Diamond_S$.
A reasonable program in this context would be a systematic comparison of diamond and tiltan, by listing combinatorial upshots of diamond and checking whether tiltan is consistent with the failure of these statements.
A classical example is the continuum hypothesis which follows from diamond.
Shelah proved in \cite{MR1623206} that tiltan is consistent with $2^\omega>\omega_1$.

Let us try to understand the difference between these combinatorial principles.
It seems that there are three basic features in which tiltan differs from diamond:
\begin{enumerate}
\item [$(\aleph)$] The prediction of tiltan is based on inclusion, while the prediction of diamond is based on equality.
\item [$(\beth)$] Diamond applies to every subset of its domain, while tiltan guesses only unbounded subsets.
\item [$(\gimel)$] A sequence of guesses of diamond is coherent, while coherency need not exist in a sequence of tiltan guesses.
\end{enumerate}
It is helpful to analyze various statements in the light of these features.
Let us consider a few examples.
If $(A_\alpha\mid\alpha\in\omega_1)$ is a diamond sequence and $A\subseteq\omega$ then $A$ is bounded in $\omega_1$.
Hence if $\delta>\omega$ satisfies $A\cap\delta=A_\delta$ then $A=A_\delta$.
It follows that every element of $[\omega]^\omega$ appears in the diamond sequence $(A_\alpha\mid\alpha\in\omega_1)$, thus $2^\omega=\omega_1$.
The consistency of tiltan with $2^\omega>\omega_1$ shows that the crucial difference here is $(\beth)$, that is the fact that bounded subsets are not predicted by tiltan.

Let us consider another example.
Galvin proved that if $2^\omega=\omega_1$ then every family $\mathcal{C}=\{C_\alpha\mid\alpha\in\omega_2\}$ of clubs of $\omega_1$ contains a subfamily $\{C_{\alpha_i}\mid i\in\omega_1\}$ so that $\bigcap_{i\in\omega_1}C_{\alpha_i}$ is a club of $\omega_1$, see \cite{MR0369081}.
Since $2^\omega=\omega_1$ follows from diamond we see that diamond implies Galvin's property.

It was shown in \cite{MR3787522} that tiltan is consistent with the failure of the Galvin property.
Since this property concentrates on unbounded subsets of $\omega_1$, it seems that $(\beth)$ is not the main factor which separates tiltan from diamond with respect to Galvin's property.
The point here is $(\gimel)$, namely the coherency of the guesses, as explicated in \cite{MR3787522}.
This means that if $(A_\alpha\mid\alpha\in\omega_1)$ is a diamond sequence, $A\subseteq\omega_1, \delta<\varepsilon<\omega_1$ and both $A\cap\delta=A_\delta, A\cap\varepsilon=A_\varepsilon$, then $A_\varepsilon\cap\delta=A_\delta$.
But if $(T_\delta\mid\delta\in\lim(\omega_1))$ is just a tiltan sequence and $T_\delta,T_\varepsilon\subseteq{A}$ then there is no any necessary connection between $T_\delta$ and $T_\varepsilon$. In fact, $T_\delta\cap T_\varepsilon=\varnothing$ is possible.

Our last example in this context is related to a question of Larson from \cite{MR2279658}.
The subject belongs to graph theory, and it is based on the following notation.
For a type $\tau$ let us say that $\tau\rightarrow(\tau,\text{infinite path})^2$ iff every graph $G$ of type $\tau$ admits either an edge-free subset of type $\tau$ or an infinite path.
Baumgartner and Larson investigated several types of graphs in the light of this relation, and one of them is $\omega^*\cdot\omega_1$.
They obtained the negative relation $\omega^*\cdot\omega_1\nrightarrow(\omega^*\cdot\omega_1,\text{infinite path})^2$ from the diamond principle (at $\aleph_1$), see \cite{MR1050558}.

Motivated by the (still open) problem whether the continuum hypothesis yields this negative relation, even if diamond fails, Larson asked whether tiltan implies $\omega^*\cdot\omega_1\nrightarrow(\omega^*\cdot\omega_1,\text{infinite path})^2$.
It is shown in \cite{gr} that tiltan is consistent with the positive relation $\omega^*\cdot\omega_1\rightarrow(\omega^*\cdot\omega_1,\text{infinite path})^2$.
Presumably, the gist of the matter here is $(\aleph)$, that is, the fact that tiltan predictions are based on inclusion only.

To see this, let us define the concept of a \emph{sofinal} set.\footnote{The term sofinal is an abbreviation of strongly cofinal from \cite{MR1050558}.}
Let $\delta$ be a limit ordinal in $\omega_1$.
A subset $S$ of $\delta\times\omega$ is sofinal in $\delta$ iff for unboundedly many $\alpha\in\delta$ there are infinitely many $n\in\omega$ so that $(\alpha,n)\in{S}$.
A subset $S$ of $\omega_1\times\omega$ is sofinal if the relevant $\delta$ is $\omega_1$.
A crucial component in the proof of the negative relation under diamond from \cite{MR1050558} is the prediction of sofinal subsets by a diamond sequence.
But the main point is that sofinal sets are predicted by sofinal subsets of $\delta$ for stationarily many $\delta\in\omega_1$.
This is done by the following consideration.
First, one fixes a bijection $h$ from $\omega_1$ to $\omega_1\times\omega$.
Second, one fixes a diamond sequence $(A_\alpha\mid\alpha\in\omega_1)$.
Finally, if $S$ is sofinal in $\omega_1\times\omega$ then the preimage of $S$ is a subset of $\omega_1$, and it is predicted by many $A_\alpha$s.
Applying $h$ to these $A_\alpha$s yields stationarily many instances in which $h''A_\alpha$ is sofinal in $\alpha$ and a subset of $S$.

The fact that sofinal sets are predicted by sofinal subsets is a key point within the proof.
This feature of a diamond sequence fades away when the prediction is based on tiltan.
Suppose that $S$ is sofinal and $X$ is the preimage of $S$ under $h$.
Let $(T_\delta\mid\delta\ \text{is a limit ordinal of}\ \omega_1)$ be a tiltan sequence, and let $Y$ be the union of all $T_\delta$s which are contained in $X$.
The set $h''Y$ is of size $\aleph_1$, but it need not be sofinal.
Actually, if tiltan holds and $2^\omega>\omega_1$ then there must be sofinal subsets for which the relevant $h''Y$ is not sofinal.
Thus, the fact that the prediction of tiltan is based on inclusion is the real reason for the difference between tiltan and diamond with respect to infinite paths in graphs.

These observations lead to the formulation of an intermediate predicting principle, similar to diamond in some sense but resembles tiltan in other features. The principle \emph{superclub} was defined by Primavesi in his dissertation, \cite{primavesi}.

\begin{definition}
\label{defsuperclub} Let $\kappa$ be regular and uncountable.
A superclub sequence for $\kappa$ is a sequence $(S_\alpha\mid\alpha\in\kappa)$ such that $S_\alpha\subseteq\alpha$ for each $\alpha\in\kappa$ and for every $A\in[\kappa]^\kappa$ there exists $B\in[A]^\kappa$ for which $S_B=\{\alpha\in\kappa\mid B\cap\alpha=S_\alpha\}$ is stationary in $\kappa$.
We shall say that superclub holds at $\kappa$ and we will denote it by $\clubsuit^{\Diamond}_\kappa$ iff there exists a superclub sequence at $\kappa$.
\end{definition}

As before, $\kappa$ can be replaced by a stationary subset $S$ of $\kappa$, in which case $\clubsuit^{\Diamond}_S$ means that there is a superclub sequence for $S$.
A comparison between superclub and its elder siblings shows that $\Diamond_S\Rightarrow\clubsuit^{\Diamond}_S\Rightarrow\clubsuit_S$.
Considering the difference between tiltan and diamond as depicted above, superclub is more akin to tiltan.
Indeed, the prediction of superclub applies only to unbounded subsets of $\kappa$, and if $A\in[\kappa]^\kappa$ is predicted by some $S_\alpha$ then $S_\alpha\subseteq A\cap\alpha$ and possibly $S_\alpha\neq A\cap\alpha$.
However, a sequence of superclub guesses is coherent.
Thus, if $A\in[\kappa]^\kappa, B\in[A]^\kappa$ and $S_B$ is the set of guesses, then $\gamma,\delta\in{S_B}$ with $\gamma<\delta$ implies $S_\delta\cap\gamma=(B\cap\delta)\cap\gamma=B\cap\gamma=S_\gamma$.
The point is that superclub acts like diamond once restricted to the set $B$.

Therefore, a study of the relationship between diamond and superclub on the one hand and superclub and tiltan on the other hand, can be guided by the above features.
It turns out that tiltan (at $\aleph_1$) is strictly weaker than superclub, and superclub is strictly weaker than diamond.
For example, Chen proved in \cite{MR3687435} that superclub is consistent with $2^\omega>\omega_1$.
Given the fact that superclub extrapolates only unbounded subsets of $\omega_1$, this result is expected.
Likewise, it was shown in \cite{MR3787522} that superclub implies Galvin's property at clubs of $\aleph_1$, and this result is expected as well since the guesses of superclub are coherent.
It shows, however, that superclub is strictly stronger than tiltan.

The purpose of this paper is two-fold.
Firstly, we shall separate tiltan from superclub at cardinals larger than $\aleph_1$.
This will be done by analyzing the Galvin property at normal filters over $\kappa=\cf(\kappa)>\aleph_1$.
Secondly, we shall prove that superclub at $\aleph_1$ implies $\mathfrak{s}=\omega_1$, where $\mathfrak{s}$ is the splitting number.
It is an interesting open problem whether tiltan implies $\mathfrak{s}=\omega_1$.
Hence we cannot be sure that this property separates tiltan from superclub.
Nevertheless, if $\kappa$ is supercompact then one can force tiltan at $\kappa^+$ with $\mathfrak{s}_\kappa>\kappa^+$.
Since superclub implies in this case that $\mathfrak{s}_\kappa=\kappa^+$, we obtain another separating property between tiltan and superclub for large successor cardinals.

Our notation is mostly standard.
We employ the Jerusalem forcing notation, thus we force upwards.
We mention here (a simple instance of) the polarized partition relation.
Let $\binom{\alpha}{\beta}\rightarrow\binom{\gamma}{\delta}$ express the statement that for every coloring $c:\alpha\times\beta\rightarrow{2}$ there are $A\in[\alpha]^\gamma,B\in[\beta]^\delta$ so that $c\upharpoonright(A\times{B})$ is constant.
Other notation will be introduced as the need arises.

\newpage

\section{Prediction principles and the Galvin property}

In this section we wish to separate $\clubsuit^\Diamond_{\kappa^+}$ from $\clubsuit_{\kappa^+}$.
Let $\mathscr{D}_{\kappa^+}$ be the club filter of $\kappa^+$.
Our first theorem shows that $\clubsuit^\Diamond_{\kappa^+}$ implies Galvin's property at $\mathscr{D}_{\kappa^+}$.
This assertion will be denoted by ${\rm Gal}(\mathscr{D}_{\kappa^+},\kappa^+,\kappa^{++})$.

Ahead of the proof we quote a basic lemma from the work of Primavesi.
Since we wish to apply this lemma to cardinals above $\aleph_1$, we spell-out the proof (in the original form it is phrased only at $\aleph_1$).
In the proof below we enumerate the elements of the superclub sequence by the limit ordinals of $\kappa^+$.
There is no loss of generality here since the prediction is done over a stationary set, so we may safely assume that all of its elements are limit ordinals.

\begin{lemma}
\label{lemclub} Assume $\clubsuit^\Diamond_{\kappa^+}$.
There exists a sequence $(B_\delta\mid\delta\in\lim(\kappa^+))$ such that $B_\delta\subseteq\delta$ is a club subset of $\delta$ for every $\delta\in\lim(\kappa^+)$, and for every club $C\subseteq\kappa^+$ one can find a set $D\subseteq{C}$ that is a club of $\kappa^+$ such that the set $\{\delta\in\lim(\kappa^+)\mid D\cap\delta=B_\delta\}$ is a stationary subset of $\kappa^+$.
\end{lemma}

\par\noindent\emph{Proof}. \newline
Let $(A_\delta\mid\delta\in\lim(\kappa^+))$ be a superclub sequence.
For each $\delta\in\lim(\kappa^+)$ let $B_\delta=c\ell(A_\delta)$.
Fix a club $C\subseteq\kappa^+$.
Apply $\clubsuit^\Diamond_{\kappa^+}$ to find $A\subseteq{C}, |A|=\kappa^+$ for which $\{\delta\in\lim(\kappa^+)\mid A\cap\delta=A_\delta\}$ is a stationary subset of $\kappa^+$.
Let $D=c\ell(A)$, so $D\subseteq{C}$ and $D$ is a club subset of $\kappa^+$.
Suppose that $\gamma\in\{\delta\in\lim(\kappa^+)\mid A\cap\delta=A_\delta\}$.
Notice that $A_\gamma\subseteq\gamma$ and $A_\gamma\subseteq{A}$.
Hence $B_\gamma=c\ell(A_\gamma)\subseteq c\ell(A)=D$.
Moreover, $B_\gamma\cap\gamma=c\ell(A)\cap\gamma=D\cap\gamma$.
Thus, $\{\gamma\in\lim(\kappa^+)\mid D\cap\gamma=B_\gamma\}$ is a stationary subset of $\kappa^+$, as required.

\hfill \qedref{lemclub}

Based on the above lemma, we proceed to the following:

\begin{theorem}
\label{thmgalsuperclub}
$\clubsuit^\Diamond_{\kappa^+}$ implies ${\rm Gal}(\mathscr{D}_{\kappa^+},\kappa^+,\kappa^{++})$.
\end{theorem}

\par\noindent\emph{Proof}. \newline
Let $(B_\delta\mid\delta\in\lim(\kappa^+))$ be a closed superclub sequence, as guaranteed by Lemma \ref{lemclub}.
Suppose that $\mathcal{C}=\{C_\alpha\mid\alpha\in\kappa^{++}\}\subseteq \mathscr{D}_{\kappa^+}$.
For every $\alpha\in\kappa^{++}$ choose $D_\alpha\subseteq{C_\alpha}$ such that $D_\alpha\in\mathscr{D}_{\kappa^+}$ and $S_\alpha=\{\delta\in\lim(\kappa^+)\mid D_\alpha\cap\delta=B_\delta\}$ is stationary in $\kappa^+$.

For every $\alpha\in\kappa^{++}$ and every $\delta\in\kappa^+$ let $H_{\alpha\delta}$ be the set $\{\beta\in\kappa^{++}\mid D_\beta\cap\delta=B_\delta\}$ provided that $D_\alpha\cap\delta=B_\delta$, and $\varnothing$ otherwise.
Thus, if $D_\alpha\cap\delta=B_\delta$ and $\beta\in{H_{\alpha\delta}}$ then $D_\beta\cap\delta=D_\alpha\cap\delta$.
Observe that $H_{\alpha\delta}\neq\varnothing$ for a stationary set of $\delta\in\kappa^+$, since $D_\alpha\cap\delta=B_\delta$ implies (trivially) that $\alpha\in{H_{\alpha\delta}}$, which means in particular that $H_{\alpha\delta}$ is not empty.
Furhtermore, if $\gamma<\delta$ and both $H_{\alpha\gamma},H_{\alpha\delta}$ are not empty then $H_{\alpha\gamma}\supseteq H_{\alpha\delta}$.
Indeed, if $\beta\in H_{\alpha\delta}$ then $D_\beta\cap\delta=D_\alpha\cap\delta$.
Hence $D_\beta\cap\gamma=(D_\beta\cap\delta)\cap\gamma =(D_\alpha\cap\delta)\cap\gamma=B_\gamma$, so $\beta\in H_{\alpha\gamma}$.

The above observations mean that every $\alpha\in\kappa^{++}$ induces a sequence $(H_{\alpha\delta}\mid\delta\in\kappa^+)$ which behaves like a tower over a stationary set of $\delta\in\kappa^+$.
The main point now is that for some $\alpha\in\kappa^{++}$ this tower consists of large sets, that is, $|H_{\alpha\delta}|=\kappa^{++}$ for stationarily many $\delta\in\kappa^+$.

For suppose not.
For every $\alpha\in\kappa^{++}$ let $\delta(\alpha)\in\kappa^+$ be the first ordinal for which $H_{\alpha\delta(\alpha)}\neq\varnothing$ and $|H_{\alpha\delta(\alpha)}|<\kappa^{++}$.
By the pigeonhole principle there is a set $A\in[\kappa^{++}]^{\kappa^{++}}$ and a fixed ordinal $\delta\in\kappa^+$ such that $\alpha\in{A}\Rightarrow\delta(\alpha)=\delta$.
Thus, $|H_{\alpha\delta}|<\kappa^{++}$ for every $\alpha\in{A}$, but also $H_{\alpha\delta}\neq\varnothing$.
Since $H_{\alpha\delta}\neq\varnothing$, it does not depend on $\alpha$ but only on $B_\delta$.
Since the number of $B_\delta$s is $\kappa^+$ and $|A|=\kappa^{++}$, there is a fixed $B$ and some $T\in[A]^{\kappa^{++}}$ such that $\alpha\in{T}\Rightarrow H_{\alpha\delta}=B$.
Consequently, $\bigcup_{\alpha\in{T}}H_{\alpha\delta}=B$ and, in particular, $|\bigcup_{\alpha\in{T}}H_{\alpha\delta}|<\kappa^{++}$.
On the other hand, $\alpha\in H_{\alpha\delta}$ for every $\alpha\in{T}$, hence $T\subseteq\bigcup_{\alpha\in{T}}H_{\alpha\delta}$.
Since $|T|=\kappa^{++}$ one concludes that $|\bigcup_{\alpha\in{T}}H_{\alpha\delta}|=\kappa^{++}$, a contradiction.

Fix, therefore, an ordinal $\alpha\in\kappa^{++}$ such that $|H_{\alpha\delta}|=\kappa^{++}$ whenever $H_{\alpha\delta}\neq\varnothing$.
Enumerate the non-empty $H_{\alpha\delta}$s by $\{H_{\alpha\varepsilon}\mid\varepsilon\in\kappa^+\}$.
By induction on $\varepsilon\in\kappa^+$ choose $\beta_\varepsilon\in H_{\alpha\varepsilon+1}$ such that $\beta_\varepsilon\notin\{\beta_\zeta\mid\zeta\in\varepsilon\}$.
This is possible since $|H_{\alpha\varepsilon+1}|=\kappa^{++}$.

We claim that $\bigcap_{\varepsilon\in\kappa^+}D_{\beta_\varepsilon} = \Delta_{\varepsilon\in\kappa^+}D_{\beta_\varepsilon}\cap D_\alpha$.
This is sufficient since $\bigcap_{\varepsilon\in\kappa^+}C_{\beta_\varepsilon} \supseteq \bigcap_{\varepsilon\in\kappa^+}D_{\beta_\varepsilon}$.
So let us prove the above equality by double inclusion.
For the easy direction, suppose that $\gamma\in\bigcap_{\varepsilon\in\kappa^+}D_{\beta_\varepsilon}$.
Certainly, $\gamma\in\Delta_{\varepsilon\in\kappa^+}D_{\beta_\varepsilon}$.
Pick $\varepsilon>\gamma$.
Notice that $\gamma\in D_{\beta_\varepsilon}\cap(\varepsilon+1)$.
But $D_{\beta_\varepsilon}\cap(\varepsilon+1)=B_{\delta(\varepsilon+1)} = D_\alpha\cap(\varepsilon+1)$, where $\delta(\varepsilon+1)$ corresponds to $\varepsilon+1$ in the enumeration of the $H_{\alpha\delta}$s which are not empty.
Thus $\gamma\in D_\alpha$, as sought.

Assume now that $\gamma\in
\Delta_{\varepsilon\in\kappa^+}D_{\beta_\varepsilon}\cap D_\alpha$.
We must show that $\gamma\in D_{\beta_\varepsilon}$ for every $\varepsilon\in\kappa^+$.
Now if $\varepsilon<\gamma$ then $\gamma\in D_{\beta_\varepsilon}$ by the definition of diagonal intersection.
If $\varepsilon\geq\gamma$ then $\gamma<\varepsilon+1$, so $\gamma\in D_\alpha\cap(\varepsilon+1)=B_{\delta(\varepsilon+1)} = D_{\beta_\varepsilon}\cap(\varepsilon+1)$, hence once again $\gamma\in D_{\beta_\varepsilon}$.
This completes the proof.

\hfill \qedref{thmgalsuperclub}

In order to separate $\clubsuit^\Diamond_{\kappa^+}$ from $\clubsuit_{\kappa^+}$, we shall prove that $\clubsuit_{\kappa^+}$ is consistent with the failure of the Galvin property at $\kappa^+$.
We consider two cases.
The simpler is the case of a successor of a regular cardinal.
In the other case, a successor of a singular cardinal, we need a bit more.
In particular, we will make use of the following lemma.
We are grateful to Yair Hayut for a very helpful discussion concerning the proof of the lemma.

\begin{lemma}
\label{lemprikry} Let $\kappa$ be a measurable cardinal and suppose that $\clubsuit_S$ holds where $S=S^{\kappa^+}_\omega$.
Let $\mathbb{P}$ be Prikry forcing at $\kappa$ and let $G\subseteq\mathbb{P}$ be generic over the ground model.
Then $\clubsuit_S$ holds in $V[G]$.
\end{lemma}

\par\noindent\emph{Proof}. \newline
Fix a tiltan sequence $(T_\delta\mid\delta\in{S})$ in $V$.
We may assume that $otp(T_\delta)=\omega$ for every $\delta\in{S}$.
Let $A$ be an unbounded subset of $\kappa^+$ in $V[G]$, and let $\name{A}$ be a name of $A$.
Fix a condition $p=(s^p,A^p)\in\mathbb{P}$ so that $p\Vdash\name{A}\in[\kappa^+]^{\kappa^+}$.

For every $\alpha\in\kappa^+$ let $D_\alpha=\{q\in\mathbb{P}\mid p\leq{q}$ and $q\Vdash\beta_\alpha=\min(\name{A}-\alpha)\}$.
Each $D_\alpha$ is a dense open subset of $\mathbb{P}$.
By the strong Prikry property, for each $\alpha\in\kappa^+$ there is a pure extension $p\leq^*p_\alpha=(s^p,A^{p_\alpha})$ and there is $n_\alpha\in\omega$ such that for every sequence $(\nu_0,\ldots,\nu_{n_\alpha-1})\in A^{p_\alpha}$ one has $p_\alpha^\frown(\nu_0,\ldots,\nu_{n_\alpha-1})\Vdash \beta_\alpha=\min(\name{A}-\alpha)$.
For every $\alpha\in\kappa^+$ we have $n_\alpha\in\omega$, hence without loss of generality $n_\alpha=n$ for each $\alpha\in\kappa^+$ and some fixed $n\in\omega$.

Let $B=\{\beta_\alpha\mid\alpha\in\kappa^+\}$.
For every $\bar{\nu}\in[\kappa]^n$ let $B_{\bar{\nu}}$ be the set $\{\beta\in\kappa^+\mid \exists q, p^\frown(\bar{\nu})\leq^* q, q\Vdash\beta\in{B}\}$.
Since $|[\kappa]^n|=\kappa<\kappa^+$, there is some $\bar{\nu}\in[\kappa]^n$ such taht $|B_{\bar{\nu}}|=\kappa^+$.
Since $B_{\bar{\nu}}\in{V}$, there is $\delta\in{S}$ for which $T_\delta\subseteq B_{\bar{\nu}}$.
Suppose that $T_\delta=\{\gamma_n\mid n\in\omega\}$.
For every $n\in\omega$ let $q_n$ be such that $p\leq^* q_n\Vdash\check{\gamma}_n\in\name{A}$.
Denote each $q_n$ by $(s^p,A^{q_n})$, and let $A=\bigcap_{n\in\omega}A^{q_n}$.
Verify that $q\in\mathbb{P}$ and $p\leq^* q\Vdash T_\delta\subseteq\name{A}$, so we are done.

\hfill \qedref{lemprikry}

There is an alternative proof of the preservation of tiltan under Prikry extensions.
Let $\mathscr{U}$ be a normal ultrafilter over $\kappa$, and let $(M_i\mid i\in\omega)$ be the associated sequence of iterates, starting from the embedding $\jmath$ derived from $\mathscr{U}$.
Let $M_\omega$ be the direct limit of the sequence and the embeddings generated from $\mathscr{U}$.
Finally, let $\bar{\kappa}=(\kappa_i\mid i\in\omega)$ be the sequence of critical points.

It was discovered by Bukowsk\'y in \cite{MR446978} and, independently, by Dehornoy in \cite{MR514228} that $\bigcap_{i\in\omega}M_i=M_\omega[\bar{\kappa}]$.
Many statements concerning Prikry forcing can be proved using the above equality, see e.g. \cite{MR4231611}.
The preservation of tiltan can serve as another example.
Suppose that there is a tiltan sequence over $S$ in the ground model, where $S=S^{\kappa^+}_\omega$ and $\kappa$ is a measurable cardinal.
If $A\in[\kappa^+]^{\kappa^+}$ is a new set in $M_\omega[\bar{\kappa}]$ then each element of $A$ comes from some $M_i$.
Hence for some $B\in[A]^{\kappa^+}$ there will be $i\in\omega$ so that $B\in{M_i}$.
An appropriate embedding of the tiltan sequence at $S$ will guess $B$ and hence $A$ as required.

We need an additional lemma about the preservation of tiltan, this time under $\aleph_1$-complete forcing notions.
It seems that this lemma is interesting in its own right, and it will also come in handy in the proof of our next theorem.
We mention here the preparatory forcing of Laver from \cite{MR0472529} to make a supercompact cardinal indestructible under $\kappa$-directed-closed forcing notions.

\begin{lemma}
\label{lemimmunity} Let $\kappa$ be an infinite cardinal.
One can obtain $\clubsuit_S$ for the set $S=S^{\kappa^{++}}_\omega$ so that every further extension of the universe by any $\aleph_1$-complete forcing notion which preserves $\kappa^{++}$ will preserve $\clubsuit_S$.
Moreover, one can get this setting where $\kappa$ is a Laver-indestructible supercompact cardinal, while preserving supercompactness and indestructibility.
\end{lemma}

\par\noindent\emph{Proof}. \newline
The statement is proved in \cite[I, \S 7]{MR1623206} and in \cite{MR3787522} in the case of $\kappa=\aleph_0$, that is $S=S^{\aleph_2}_{\aleph_0}$.
The general case is identical, and since this forcing is $\kappa$-directed-closed one obtains a preservation of supercompactness and Laver-indestructibility.

\hfill \qedref{lemimmunity}

We can prove now the consistency of tiltan with the failure of the Galvin property at arbitrarily large successor cardinals.

\begin{theorem}
\label{thmtiltanegalvin} Let $\kappa$ be a regular cardinal.
Then one can force $\clubsuit_{\kappa^+}$ with $\neg{\rm Gal}(\mathscr{D}_{\kappa^+},\kappa^+,\kappa^{++})$.
If there exists a supercompact cardinal then one can force $\kappa$ to be singular, with $\clubsuit_{\kappa^+}$ and $\neg{\rm Gal}(\mathscr{D}_{\kappa^+},\kappa^+,\kappa^{++})$.
\end{theorem}

\par\noindent\emph{Proof}. \newline
Fix a regular cardinal $\kappa$ and suppose that $\kappa^{++}$ carries an indestructible tiltan sequence at $S=S^{\kappa^{++}}_\omega$ with respect to $\aleph_1$-complete forcing notions.
Such a sequence can be forced by virtue of Lemma \ref{lemimmunity}.
We may assume that $2^\kappa=\kappa^+$ and $2^{\kappa^+}=\kappa^{++}$.

Apply the forcing of Abraham and Shelah from \cite{MR830084} to $\kappa^{++}$, and call the generic extension $V[G]$.
Thus in $V[G]$ there is some $\lambda>\kappa^{++}$ and a family $\mathcal{C}=\{C_\alpha\mid\alpha\in\lambda\}\subseteq\mathscr{D}_{\kappa^{++}}$ such that $\mathcal{C}$ witnesses $\neg{\rm Gal}(\mathscr{D}_{\kappa^{++}},\kappa^{++},\kappa^{+3})$.
Since the forcing of Abraham and Shelah is $\aleph_1$-complete (and much more) we see that $\clubsuit_S$ holds in $V[G]$.

Let $\mathbb{Q}$ be the collapse of $\kappa^+$ to $\kappa$ as defined in $V[G]$ and let $H\subseteq\mathbb{Q}$ be generic over $V[G]$.
Denote $S^{\kappa^+}_\omega$ as computed in $V[G]$ by $T$.
Observe that both $\neg{\rm Gal}(\mathscr{D}_{\kappa^+},\kappa^+,\kappa^{++})$ and $\clubsuit_T$ are true in $V[G][H]$, as proved in \cite[Theorem 2.7]{MR3787522}.
This completes the proof of the first statement of the theorem as $\clubsuit_T$ implies $\clubsuit_{\kappa^+}$.

Assume now that $\kappa$ is supercompact in the ground model.
Force $\kappa$ to be Laver-indestructible.
Apply Lemma \ref{lemimmunity} to obtain a tiltan sequence at $S=S^{\kappa^{++}}_\omega$ which is indestructible under $\aleph_1$-complete forcing extensions.
As a second step, force $\neg{\rm Gal}(\mathscr{D}_{\kappa^+},\kappa^+,\kappa^{++})$ with $\clubsuit_T$ at $T=S^{\kappa^+}_\omega$ as done in the first part of the proof.
Notice that $\kappa$ is supercompact in $V[G][H]$ since both the Abraham-Shelah forcing and the collapse of $\kappa^+$ to $\kappa$ are $\kappa$-directed-closed.

Let $\mathbb{P}$ be Prikry forcing into $\kappa$ as defined in $V[G][H]$, and let $K\subseteq\mathbb{P}$ be $V[G][H]$-generic.
From \cite{bgp} we know that $\neg{\rm Gal}(\mathscr{D}_{\kappa^+},\kappa^+,\kappa^{++})$ holds in $V[G][H][K]$. From Lemma \ref{lemprikry} we infer that $\clubsuit_T$ holds in $V[G][H][K]$.
Thus the proof is accomplished.

\hfill \qedref{thmtiltanegalvin}

We can phrase now the main result of this section:

\begin{corollary}
\label{corseparation} The principles $\clubsuit_{\kappa^+}$ and $\clubsuit^\Diamond_{\kappa^+}$ are separable for every $\kappa=\cf(\kappa)$, and consistently they are separable at successors of singular cardinals.
\end{corollary}

It is plausible that $\clubsuit_{\kappa^+}$ and $\clubsuit^\Diamond_{\kappa^+}$ are separable in a global sense.
This will be especially interesting at successors of singular cardinals, as dictated in the following problem.

\begin{question}
\label{qglobal} Is it consistent that $\clubsuit_{\lambda^+}$ holds and $\clubsuit^\Diamond_{\lambda^+}$ fails at every singular cardinal $\lambda$ simultaneously?
\end{question}

\newpage

\section{Superclub and the splitting number}

In this section we show that superclub implies $\mathfrak{s}=\omega_1$.
More generally, if $\kappa$ is weakly compact then the statement $\mathfrak{s}_\kappa=\kappa^+$ follows from $\clubsuit^\Diamond_{\kappa^+}$.
Let us commence with the pertinent definition.

\begin{definition}
\label{defsplitting} The splitting number. \newline
Let $\kappa$ be an infinite cardinal.
\begin{enumerate}
\item [$(\aleph)$] For $B,S\in[\kappa]^\kappa$ one says that $S$ splits $B$ iff $|B\cap{S}|=|B\cap(\kappa-S)|=\kappa$.
\item [$(\beth)$] A family $\{S_\alpha\mid\alpha\in\lambda\}\subseteq[\kappa]^\kappa$ is a splitting family iff every $B\in[\kappa]^\kappa$ is split by some $S_\alpha$.
\item [$(\gimel)$] The splitting number $\mathfrak{s}_\kappa$ is the minimal cardinality of a splitting family for $[\kappa]^\kappa$.
\end{enumerate}
\end{definition}

As usual, if $\kappa=\aleph_0$ then $\mathfrak{s}_\kappa$ is denoted by $\mathfrak{s}$.
The splitting number belongs to the family of cardinal characteristics.
All of them stabilize at $\omega_1$ under the continuum hypothesis, a fortiori under diamond.
Hence it is a good test problem to check the possible value of any cardinal characteristic under tiltan.

For small characteristics the value is $\omega_1$, e.g. tiltan implies $\mathfrak{p}=\omega_1$.
For large characteristics, tiltan is consistent with a large value, e.g. one can force tiltan with $\mathfrak{d}>\omega_1$.
Medium characteristics are more interesting and, in particular, it is open whether tiltan is consistent with $\mathfrak{s}>\omega_1$.
Our primary goal in this section is to prove that superclub implies $\mathfrak{s}=\omega_1$.
Let us indicate that for some characteristics, a large value is consistent with superclub, see \cite{1179}.
In this light, the small value of the splitting number becomes more interesting.

Inasmuch as we do not know whether $\mathfrak{s}>\omega_1$ is consistent with tiltan, we cannot separate tiltan at $\aleph_1$ from superclub at $\aleph_1$ on the grounds of the splitting number.
However, we can do this at larger cardinals.
It was shown by Suzuki in \cite{MR1251349} that $\mathfrak{s}_\kappa>\kappa$ iff $\kappa$ is weakly compact.
We shall see that $\clubsuit^\Diamond_{\kappa^+}$ implies $\mathfrak{s}_\kappa=\kappa^+$ for weakly compact cardinals.
On the other hand, $\clubsuit_{\kappa^+}$ is consistent with $\mathfrak{s}_\kappa>\kappa^+$ for sufficiently large cardinals, as proved in \cite[Theorem 4.5]{MR3787522}.
Thus another way to separate tiltan from superclub (at successors of large cardinals) is obtained.

\begin{theorem}
\label{thmsuperclubsplitting} Superclub implies $\mathfrak{s}=\omega_1$.
\end{theorem}

\par\noindent\emph{Proof}. \newline
From \cite[Claim 2.4]{MR2927607} we know that if $\mathfrak{s}>\omega_1$ then $\binom{\omega_1}{\omega}\rightarrow\binom{\omega_1}{\omega}$.
We shall prove that superclub implies the negative relation $\binom{\omega_1}{\omega}\nrightarrow\binom{\omega_1}{\omega}$ and conclude that $\mathfrak{s}=\omega_1$.
Let $(S_\delta\mid\delta\in\omega_1)$ be a superclub sequence.
For every $\ell\in\{0,1\}$ let $S^\ell=\{\delta\in\omega_1\mid(\forall\alpha\in{S_\delta})(\alpha=\ell\ {\rm mod}\ 2)\}$.
Notice that $S^\ell\neq\varnothing$ for $\ell\in\{0,1\}$, since the set of even ordinals is predicted by the superclub sequence, as well as the set of odd ordinals.

By induction on $\alpha\in\omega_1$ we choose a set $A_\alpha\subseteq\omega$ such that:
\begin{enumerate}
\item [$(a)$] $|A_\alpha|=|\omega-A_\alpha|=\aleph_0$.
\item [$(b)$] If $\delta\in\alpha\cap{S^0}$ and $B^0_\delta=\bigcap\{(\omega-A_\beta)\mid 2\beta\in S_\delta\}$ is infinite, then there are $n^0_\delta,n^1_\delta\in B^0_\delta$ such that $n^0_\delta\notin A_\alpha, n^1_\delta\in A_\alpha$.
\item [$(c)$] If $\delta\in\alpha\cap{S^1}$ and $B^1_\delta=\bigcap\{A_\beta\mid 2\beta+1\in S_\delta\}$ is infinite, then there are $n^0_\delta,n^1_\delta\in B^1_\delta$ such that $n^0_\delta\notin A_\alpha, n^1_\delta\in A_\alpha$.
\end{enumerate}
The choice is possible since for every $\alpha\in\omega_1$ one has $|\alpha\cap S^0|,|\alpha\cap S^1|\leq\aleph_0$.
Thus the number of $B^\ell_\delta$s to consider is at most $\aleph_0$ and by induction on $\omega$ one can pick two new distinct elements from $B^\ell_\delta$ and send them to $A_\alpha$ and $\omega-A_\alpha$ respectively.

Define $c:\omega_1\times\omega\rightarrow{2}$ by letting $c(\alpha,n)=1$ iff $n\in{A_\alpha}$.
Assume towards contradiction that there are $Y\in[\omega_1]^{\omega_1}, Z\in[\omega]^\omega$ and $\ell\in\{0,1\}$ such that $c''(Y\times{Z})=\{\ell\}$.
Let $Y=\{\alpha_j\mid j\in\omega_1\}$, let $Z$ be maximal with respect to the property that $c''(Y\times{Z})=\{\ell\}$ and let $X=\{2\alpha+\ell\mid\alpha\in{Y}\}$.
Let $W\in[X]^{\aleph_1}$ be a set on which the superclub sequence acts like diamond.
Set $T=\{\delta\in\omega_1\mid W\cap\delta=S_\delta\}$, so $T$ is a stationary subset of $\omega_1$.
For each $\alpha\in\omega_1$ let $A^1_\alpha=A_\alpha$ and $A^0_\alpha=\omega-A_\alpha$.
Choose a sufficiently large $\delta\in{T}$ such that:
\[\bigcap\{A^\ell_\alpha\mid\alpha\in{Y}\} = \bigcap\{A^\ell_\alpha\mid\alpha\in{Y}\cap\delta\}\]
There is such $\delta$ since the sequence of sets $(\bigcap\limits_{j<i}A^\ell_{\alpha_j}\mid i\in\omega_1)$ is decreasing with respect to inclusion, and being a sequence of subsets of $\omega$ it stabilizes after countably many steps.

Notice that $\bigcap\{A^\ell_\alpha\mid\alpha\in{Y}\cap\delta\} = B^\ell_\delta$ since $W\cap\delta=S_\delta$.
We know that this set is infinite since it contains $Z$, and we may assume that it equals $Z$ by taking $Z$ to be maximal with respect to inclusion.
Fix $\alpha\in{Y}-(\delta+1)$.
By the construction of $(A_\alpha\mid\alpha\in\omega_1)$ there are $n^0_\delta,n^1_\delta\in B^\ell_\delta=Z$ such that $n^0_\delta\notin A_\alpha$ and $n^1_\delta\in A_\alpha$.
By definition, $c(\alpha,n^0_\delta)=0$ and $c(\alpha,n^1_\delta)=1$, thus $c''(Y\times{Z})=\{0,1\}$, a contradiction.
Since $\mathfrak{s}>\omega_1$ implies $\binom{\omega_1}{\omega}\rightarrow\binom{\omega_1}{\omega}$ we deduce that superclub implies $\mathfrak{s}=\omega_1$, as sought.

\hfill \qedref{thmsuperclubsplitting}

The above proof generalizes to higher cardinals.
Of course, it would be interesting only if $\mathfrak{s}_\kappa>\kappa$, namely in cases of weak compactness.

\begin{claim}
\label{clmgeneralized} Let $\kappa$ be a weakly compact cardinal.
Then $\clubsuit^\Diamond_{\kappa^+}$ implies $\mathfrak{s}_\kappa=\kappa^+$.
\end{claim}

\par\noindent\emph{Proof}. \newline
From a straightforward generalization of the main theorem of this section we infer that $\clubsuit^\Diamond_{\kappa^+}$ implies the negative relation $\binom{\kappa^+}{\kappa}\nrightarrow\binom{\kappa^+}{\kappa}$ where $\kappa$ is weakly compact.
From \cite[Claim 1.2]{MR3201820} we know that if $\mathfrak{s}_\kappa>\kappa^+$ then $\binom{\kappa^+}{\kappa}\rightarrow\binom{\kappa^+}{\kappa}$.
Thus under superclub at $\kappa^+$ necessarily $\mathfrak{s}_\kappa=\kappa^+$, as required.

\hfill \qedref{clmgeneralized}

We can derive now the following conclusion which gives an alternative way to separate tiltan from superclub.

\begin{corollary}
\label{corsplitting} Let $\kappa$ be supercompact.
Then one can force $\clubsuit_{\kappa^+}$ with $\neg\clubsuit^\Diamond_{\kappa^+}$.
\end{corollary}

\par\noindent\emph{Proof}. \newline
It is proved in \cite[Theorem 4.5]{MR3787522} that $\clubsuit_{\kappa^+}$ is consistent with $\mathfrak{s}_\kappa>\kappa^+$ where $\kappa$ is supercompact.
By Claim \ref{clmgeneralized}, $\clubsuit^\Diamond_{\kappa^+}$ fails in this model, so we are done.

\hfill \qedref{corsplitting}

One may wonder whether Galvin's property and the splitting number are the only statements which separate tiltan from superclub.
Formally, suppose that $\clubsuit_{\kappa^+}$ holds and $\clubsuit^\Diamond_{\kappa^+}$ fails.
Does it follow that either $\mathfrak{s}_\kappa>\kappa^+$ or $\neg{\rm Gal}(\mathscr{D}_{\kappa^+},\kappa^+,\kappa^{++})$?
We believe that the answer is negative.

A possible strategy towards a negative answer would be as follows.
Begin with a supercompact cardinal $\kappa$, force $\mathfrak{s}_\kappa>\kappa^+$ but keep the Galvin property.
This is possible, at least for some instances of the Galvin proeprty.
Moreover, $\clubsuit_{\kappa^+}$ is forceable with this setting.
In such models $\clubsuit^\Diamond_{\kappa^+}$ necessarily fails, since it implies $\mathfrak{s}_\kappa=\kappa^+$, thus $\clubsuit^\Diamond_{\kappa^+}$ is separated from $\clubsuit_{\kappa^+}$.
Now force Prikry into $\kappa$.
By \cite{bgp}, Galvin's property holds in the generic extension.
By Lemma \ref{lemprikry}, $\clubsuit_{\kappa^+}$ holds as well.
Since $\kappa$ remains a strong limit cardinal and becomes singular, $\mathfrak{s}_\kappa=\mathfrak{s}_{\cf(\kappa)}<\kappa$.
The missing part is to show that Prikry forcing does not add a superclub sequence.
This is plausible, but it requires a formal argument.

\newpage

\section{Open problems}

In this section we introduce some open problems which arise naturally from the results in the previous sections.
The first problem has a local flavor, but it is relatively old and it seems quite stubborn.

\begin{question}
\label{qtiltansplitting} Is it consistent that tiltan holds and $\mathfrak{s}>\omega_1$?
\end{question}

This problem is strongly related to the question of possible consistency of tiltan with the polarized relation $\binom{\omega_1}{\omega}\rightarrow\binom{\omega_1}{\omega}$.
We indicate that $\mathfrak{s}>\omega_1$ implies $\binom{\omega_1}{\omega}\rightarrow\binom{\omega_1}{\omega}$, but not vice versa, so these statements are not equivalent.

Of course, diamond implies $\mathfrak{s}=\omega_1$ since diamond implies $2^\omega=\omega_1$.
However, $\mathfrak{s}=\omega_1$ is a weaker statement than $2^\omega=\omega_1$.
Although $\mathfrak{s}=\omega_1$ is an assertion about sets of natural numbers, superclub implies $\mathfrak{s}=\omega_1$.
Thus it seems that the prediction of countable sets by a diamond sequence is not the crucial factor with respect to the statement $\mathfrak{s}=\omega_1$.
We tend to believe that $\mathfrak{s}=\omega_1$ does not follow from tiltan.

There are several interesting problems around tiltan and superclub with respect to other cardinal characteristics, but we focus here on the splitting number.
Let us introduce another interesting problem in this context.

\begin{question}
\label{qtiltansplkappa} Let $\kappa$ be a supercompact cardinal and let $S=S^{\kappa^+}_\kappa$.
Is it consistent that $\clubsuit_S$ holds and $\mathfrak{s}_\kappa>\kappa^+$?
\end{question}

We indicate that the consistency of $\clubsuit_{\kappa^+}$ with $\mathfrak{s}_\kappa>\kappa^+$ as proved in \cite{MR3787522} comes from $\clubsuit_S$ where $S=S^{\kappa^+}_\theta$ for some $\theta\in\kappa$.
The difficulties in the case of $S=S^{\kappa^+}_\kappa$ are similar to the difficulty of tiltan at $\aleph_1$ with $\mathfrak{s}>\omega_1$.

Our next topic is Luzin sequences and inseparable families.
Suppose that $\mathcal{A}=\{A_\beta:\beta\in\omega_1\}\subseteq[\omega]^\omega$.
A set $B\in[\omega]^\omega$ \emph{separates} $\mathcal{A}$ iff $|\{\beta\in\omega_1:A_\beta\subseteq^* B\}|= |\{\beta\in\omega_1:A_\beta\subseteq^*(\omega-B)\}|=\aleph_1$.
A family $\mathcal{A}\subseteq[\omega]^\omega$ will be called \emph{inseparable} iff $\mathcal{A}$ is an almost disjoint family and no $B\in[\omega]^\omega$ separates $\mathcal{A}$.
A family $\mathcal{A}=\{A_\beta:\beta\in\omega_1\}\subseteq[\omega]^\omega$ is a \emph{Luzin sequence} iff $\mathcal{A}$ is almost disjoint and for every $\alpha\in\omega_1,n\in\omega$ the set $\{\beta\in\alpha:A_\alpha\cap A_\beta\subseteq n\}$ is finite.

It is evident that if $\mathcal{A}$ is a Luzin sequence then $\mathcal{A}$ is inseparable.
Luzin, \cite{MR0010601}, proved that Luzin sequences exist, in ZFC.
Hence inseparable families exist in ZFC as well.
Abraham and Shelah asked in \cite{MR1852375} whether the only way to produce inseparable families is through Luzin sequences.
They proved that the answer is independent of ZFC.
If $2^\omega=\omega_1$ then there exists an inseparable family which contains no Luzin subsequence.
On the other hand, Martin's axiom and $2^\omega>\omega_1$ implies that every inseparable family contains a Luzin subsequence.
Remark that the existence of an inseparable family which contains no Luzin subsequence is consistent with $2^\omega>\omega_1$.
One can begin with $2^\omega=\omega_1$ and add $\lambda$-many Cohen reals.
The inseparable family which contains no Luzin subsequence will not be destroyed by these Cohen reals.

\begin{question}
\label{qluzin} Tiltan principles and Inseparable families.
\begin{enumerate}
\item [$(\aleph)$] Does tiltan imply the existence of an inseparable family which contains no Luzin subsequence?
\item [$(\beth)$] Does superclub imply the existence of an inseparable family which contains no Luzin subsequence?
\end{enumerate}
\end{question}

There are two elements in the proof under $2^\omega=\omega_1$ which require attention.
The first one is an enumeration of countable subsets of $\omega_1$ which are supposed to approximate some Luzin subsequence of size $\aleph_1$.
Under $2^\omega=\omega_1$ one enumerates all of them, but this enumeration can be replaced by an appropriate tiltan sequence.
The second is an enumeration of $[\omega]^\omega$ of length $\omega_1$, needed for taking care of every $B\in[\omega]^\omega$ through an inductive process of length $\omega_1$, where at each step we have only $\aleph_0$-many sets to take care of.
This element seems problematic, as tiltan (or superclub) say nothing about countable sets.

Finally, in this work we separate tiltan from superclub at arbitrarily large cardinals.
One may wonder about a possible separation of superclub from diamond in a similar fashion.
We strongly believe that superclub at $\kappa^+$ is weaker than diamond at $\kappa^+$ where $\kappa$ is arbitrarily large.
It is known that tiltan with $2^\omega=\omega_1$ implies diamond, and a fortiori for superclub.
The same is true in general, and for the sake of completeness we give here the proof (for the case of $\aleph_1$ see, e.g., \cite[I, Fact 7.3]{MR1623206}).

\begin{claim}
\label{clmtiltandiamond} Let $S$ be a stationary subset of $\kappa^+$ which consists of limit ordinals.
Assume that $2^\kappa=\kappa^+$.
Then $\clubsuit_S$ implies $\Diamond_S$.
\end{claim}

\par\noindent\emph{Proof}. \newline
Let $(T_\alpha\mid\alpha\in{S})$ be a $\clubsuit_S$-sequence.
We elicit another sequence of sets $(D_\alpha\mid\alpha\in{S})$ and we prove that this is a diamond sequence.
Let $(B_i\mid i\in\kappa^+)$ be an enumeration of all the bounded subsets of $\kappa^+$, where each bounded subset appears $\kappa^+$ many times in the enumeration and $B_i\subseteq{i}$ for every $i\in\kappa^+$.
For every $\alpha\in{S}$ we define $D_\alpha=\bigcup_{i\in{A_\alpha}}B_i$.
We shall show that this is a diamond sequence.

Firstly, let $X$ be a bounded subset of $\kappa^+$.
By the properties of our enumeration, the set $I=\{i\in\kappa^+\mid X=B_i\}$ is unbounded in $\kappa^+$.
Hence there is a stationary subset $S_0$ of $S$ such that $T_\alpha\subseteq{I}$ for every $\alpha\in{S_0}$.
For each $\alpha\in{S_0}$ one has $D_\alpha=\bigcup_{i\in{A_\alpha}}B_i=\bigcup_{i\in{A_\alpha}}X=X$.
Pick any $\alpha>\sup(X)$ from $S_0$.
One can see that $D_\alpha=X=X\cap\alpha$, thus for an end-segment of $S_0$ one has $X\cap\alpha=D_\alpha$, so we are done.

Secondly, let $X$ be an unbounded subset of $\kappa^+$.
We define $f:\kappa^+\rightarrow\kappa^+$ by induction on $\alpha\in\kappa^+$ as follows.
Assuming that $f(\beta)$ is already defined for every $\beta\in\alpha$,
let $f(\alpha)$ be the first ordinal $i\in\kappa^+$ greater than every element of $\{f(\beta)\mid\beta\in\alpha\}$ such that $B_i=X\cap\sup(\{f(\beta)\mid\beta\in\alpha\})$.
By induction one can verify that $f(\alpha)\geq\alpha$ for every $\alpha\in\kappa^+$.
Let $C$ be a club of $\kappa^+$, closed under $f$, that is if $\delta\in{C}$ then $f(\beta)\in\delta$ for every $\beta\in\delta$.
Let $X'=\{f(\alpha)\mid\alpha\in\kappa^+\}$.
Since $f$ is increasing, $X'$ is unbounded in $\kappa^+$ and hence there is a stationary subset of ordinals $\delta\in S\cap{C}$ for which $T_\delta\subseteq{X'}$.
We shall show that for each such $\delta$ one has $D_\delta=X\cap\delta$, and this will accomplish the argument.

We shall prove the above equality by double inclusion.
Let us begin with the easy direction, namely $D_\delta\subseteq X\cap\delta$.
Assume, therefore, that $\gamma\in{D_\delta}$, and since $D_\delta=\bigcup_{i\in{T_\delta}}B_i$ there is some $i\in{T_\delta}$ for which $\gamma\in{B_i}$.
Recall that $T_\delta\subseteq{X'}$, and this means that $i\in{T_\delta}$ is of the form $f(\alpha)$ for some $\alpha\in\kappa^+$.
By definition, $B_i=X\cap\sup(\{f(\beta)\mid\beta\in\alpha\})$, and in particular $B_i\subseteq{X}$.
Since $\gamma\in{B_i}$ we conclude that $\gamma\in{X}$.
Likewise, $B_i\subseteq{i}=f(\alpha)<\delta$ since $f(\alpha)=i\in{T_\delta}\subseteq\delta$, so $\gamma\in X\cap\delta$ and the first inclusion is established.

Let us prove the other direction, that is $X\cap\delta\subseteq{D_\delta}$.
Assume that $\beta\in{X}\cap\delta$.
Firstly we show that there exists $i\in{T_\delta}$ so that $i>f(\beta)$.
Indeed, $\beta\in\delta$ and hence $f(\beta)<\delta$ as $\delta\in{C}$.
Likewise, $T_\delta$ is unbounded in $\delta$ and hence one can pick an ordinal $\gamma$ so that $f(\beta)<\gamma\in\delta$ and $\gamma\in{T_\delta}$.
Since $T_\delta\subseteq{X'}$, there is an ordinal $\alpha$ such that $\gamma=f(\alpha)$.
Notice that $\beta<\alpha$ since $f(\beta)<\gamma=f(\alpha)$, so choose $i$ to be $\gamma$.

Now the existence of such an ordinal $i$ implies that $i=f(\alpha)$ for some $\alpha\in\kappa^+$ since $T_\delta\subseteq{X'}$.
Observe that $\alpha>\beta$ since $f(\beta)<i=f(\alpha)$ and $f$ is an increasing function.
Recall that $B_i=X\cap\sup(\{f(\beta)\mid\beta\in\alpha\})$ and conclude that $\beta\in{B_i}\subseteq\bigcup_{i\in{A_\delta}}B_i=D_\delta$, and the proof is accomplished.

\hfill \qedref{clmtiltandiamond}

From the above claim we deduce that a separation between diamond and superclub at $\kappa^+$ boils down to a model of $\clubsuit^\Diamond_{\kappa^+}$ with $2^\kappa>\kappa^+$.
We expect a positive answer to the following question.

\begin{question}
\label{qdiamondsuperclub} Is it consistent that superclub holds at $\kappa^+$ and $2^\kappa>\kappa^+$ where $\kappa$ is arbitrarily large?
\end{question}

\newpage

\bibliographystyle{alpha}
\bibliography{arlist}

\begin{thebibliography}{BHM75}

\bibitem[AS86]{MR830084}
U.~Abraham and S.~Shelah.
\newblock On the intersection of closed unbounded sets.
\newblock {\em J. Symbolic Logic}, 51(1):180--189, 1986.

\bibitem[AS01]{MR1852375}
Uri Abraham and Saharon Shelah.
\newblock Lusin sequences under {CH} and under {M}artin's axiom.
\newblock {\em Fund. Math.}, 169(2):97--103, 2001.

\bibitem[BGP21]{bgp}
Tom Benhamou, Shimon Garti, and Alejandro Poveda.
\newblock Negating the {G}alvin property.
\newblock {\em preprint}, 2021.

\bibitem[BHM75]{MR0369081}
J.~E. Baumgartner, A.~H\c{a}j\c{n}al, and A.~Mate.
\newblock Weak saturation properties of ideals.
\newblock In {\em Infinite and finite sets ({C}olloq., {K}eszthely, 1973;
  dedicated to {P}. {E}rd\H{o}s on his 60th birthday), {V}ol. {I}}, pages
  137--158. Colloq. Math. Soc. J\'{a}nos Bolyai, Vol. 10. 1975.

\bibitem[BL90]{MR1050558}
James~E. Baumgartner and Jean~A. Larson.
\newblock A diamond example of an ordinal graph with no infinite paths.
\newblock {\em Ann. Pure Appl. Logic}, 47(1):1--10, 1990.

\bibitem[Buk77]{MR446978}
Lev Bukovsk\'{y}.
\newblock Iterated ultrapower and {P}rikry's forcing.
\newblock {\em Comment. Math. Univ. Carolinae}, 18(1):77--85, 1977.

\bibitem[Che17]{MR3687435}
William Chen.
\newblock Variations of the stick principle.
\newblock {\em Eur. J. Math.}, 3(3):650--658, 2017.

\bibitem[Deh78]{MR514228}
Patrick Dehornoy.
\newblock Iterated ultrapowers and {P}rikry forcing.
\newblock {\em Ann. Math. Logic}, 15(2):109--160, 1978.

\bibitem[Gar18]{MR3787522}
Shimon Garti.
\newblock Tiltan.
\newblock {\em C. R. Math. Acad. Sci. Paris}, 356(4):351--359, 2018.

\bibitem[Gar23]{gr}
Shimon Garti.
\newblock Tiltan and graphs with no infinite paths.
\newblock {\em Period. Math. Hungar.}, (To appear), 2023.

\bibitem[GS12]{MR2927607}
Shimon Garti and Saharon Shelah.
\newblock Strong polarized relations for the continuum.
\newblock {\em Ann. Comb.}, 16(2):271--276, 2012.

\bibitem[GS14]{MR3201820}
Shimon Garti and Saharon Shelah.
\newblock Partition calculus and cardinal invariants.
\newblock {\em J. Math. Soc. Japan}, 66(2):425--434, 2014.

\bibitem[GS23]{1179}
Shimon Garti and Saharon Shelah.
\newblock Tiltan and superclub.
\newblock {\em Accepted}, 2023.

\bibitem[HU20]{MR4231611}
Yair Hayut and Spencer Unger.
\newblock Stationary reflection.
\newblock {\em J. Symb. Log.}, 85(3):937--959, 2020.

\bibitem[Jen72]{MR309729}
R.~Bj\"{o}rn Jensen.
\newblock The fine structure of the constructible hierarchy.
\newblock {\em Ann. Math. Logic}, 4:229--308; erratum, ibid. 4 (1972), 443,
  1972.
\newblock With a section by Jack Silver.

\bibitem[Lar06]{MR2279658}
Jean~A. Larson.
\newblock Partition relations on a plain product order type.
\newblock {\em Ann. Pure Appl. Logic}, 144(1-3):117--125, 2006.

\bibitem[Lav78]{MR0472529}
Richard Laver.
\newblock Making the supercompactness of {$\kappa $} indestructible under
  {$\kappa $}-directed closed forcing.
\newblock {\em Israel J. Math.}, 29(4):385--388, 1978.

\bibitem[Lus43]{MR0010601}
N.~Lusin.
\newblock Sur les parties de la suite naturelle des nombres entiers.
\newblock {\em C. R. (Doklady) Acad. Sci. URSS (N.S.)}, 40:175--178, 1943.

\bibitem[Ost76]{MR438292}
A.~J. Ostaszewski.
\newblock On countably compact, perfectly normal spaces.
\newblock {\em J. London Math. Soc. (2)}, 14(3):505--516, 1976.

\bibitem[Pri11]{primavesi}
Alexander Primavesi.
\newblock Guessing axioms, invariance and {S}uslin trees.
\newblock {\em Ph.D. Thesis, University of East Anglia}, 2011.

\bibitem[She98]{MR1623206}
Saharon Shelah.
\newblock {\em Proper and improper forcing}.
\newblock Perspectives in Mathematical Logic. Springer-Verlag, Berlin, second
  edition, 1998.

\bibitem[Suz93]{MR1251349}
Toshio Suzuki.
\newblock On splitting numbers.
\newblock {\em urikaisekikenky usho Koky uroku}, (818):118--120, 1993.
\newblock Mathematical logic and applications '92 (Japanese) (Kyoto, 1992).

\end{thebibliography}

\end{document}